\documentclass [11pt]{amsart}
\usepackage{amssymb}
\usepackage{xypic}

\newcommand{\rpfeil}[2]{\stackrel{#2}{\verylongarrow{#1mm}}}

\newtheorem{theorem}{Theorem}[section]

\newtheorem{lemma}[theorem]{Lemma}
\newtheorem{proposition}[theorem]{Proposition}
\newtheorem{definition}[theorem]{Definition}
\newtheorem{corollary}[theorem]{Corollary}
\newtheorem{conjecture}[theorem]{Conjecture}
\newtheorem{exmp}[theorem]{Example}
\newtheorem{exmps}[theorem]{Examples}
\newtheorem{rem}[theorem]{Remark}
\newenvironment{example}{\begin{exmp}\rm}{\end{exmp}}

\newenvironment{remark}{\begin{rem}\rm}{\end{rem}\rm}
\newcommand{\prf}{{\em Proof}. }

\newcommand{\beeq}[1]{\begin{eqnarray}\label{#1}}
\newcommand{\eneq}{\end{eqnarray}}

\newcommand{\ka}{{\mathcal A}}
\newcommand{\ko}{{\mathcal O}}
\newcommand{\kp}{{\mathcal P}}

\newcommand{\IC}{{\mathbb C}}

\newcommand{\IP}{{\mathbb P}}
\newcommand{\IQ}{{\mathbb Q}}
\newcommand{\IR}{{\mathbb R}}
\newcommand{\IZ}{{\mathbb Z}}
\newcommand{\gm}{{\mathfrak M}}
\newcommand{\gn}{{\mathfrak N}}

\newcommand{\gggg}{{\mathfrak g}}
\newcommand{\ssl}{{\mathfrak s}{\mathfrak l}}
\newcommand{\sso}{{\mathfrak s}{\mathfrak o}}

\newcommand{\OO}{{\rm O}}

\newcommand{\verylongarrow}[1]{\hbox to #1{\rightarrowfill}}
\newcommand{\silent}[1]{}
\begin{document}
\title[Generalized Calabi--Yau structures, K3 Surfaces, and B-fields]
       {Generalized Calabi--Yau structures, K3 Surfaces, and B-fields}

\email{huybrech@math.jussieu.fr}
\author[Daniel Huybrechts]{Daniel Huybrechts}
\address{Institut de Math\'ematiques de Jussieu, 2 place de Jussieu, 75251 Paris Cedex 05, France}

\maketitle

This article collects a few observations concerning Hitchin's generalized Calabi--Yau structures
in dimension four. I became interested  in these while thinking about the moduli space
of K3 surfaces (with metric and $B$-field) and its relation to the moduli space
of $N=(2,2)$ SCFT.

Roughly, a generalized Calabi--Yau structure is a very special even non-degenerate complex form, which usually will
be called $\varphi$. The main examples are $\varphi=\sigma$, where $\sigma$ is the holomorphic two-form on a K3 surface,
and $\varphi=\exp(i\omega)$, where  $\omega$ is an arbitrary symplectic form. A generalized K3 surface consists of a pair
$(\varphi,\varphi')$ of generalized Calabi--Yau structures satisfying certain orthogonality conditions
which are modeled on the relation between the holomorphic two-form $\sigma$ on a K3 surface and a 
Ricci-flat K\"ahler form $\omega$.

As was explained by Aspinwall and Morrison (cf.\ \cite{AM, Huy, NW}), the moduli space $\gm_{(2,2)}$
of $N=(2,2)$ SCFT  fibers over the moduli space $\gm_{(4,4)}$ of $N=(4,4)$
SCFT. The fibre of the projection $\gm_{(2,2)}\to\gm_{(4,4)}$
is isomorphic to $S^2\times S^2$. Using the period map, the moduli space of $B$-field shifts of
hyperk\"ahler metrics $\gm_{{\rm HK}}$ can be identified with an open dense subset of 
$\gm_{(4,4)}$. 
For any chosen hyperk\"ahler metric $g\in\gm_{\rm HK}$ there is an $S^2$ worth of complex structures
making this metric a K\"ahler metric. Thus, the moduli space $\gm_{{\rm K3}}$  of $B$-field shifts
of complex K3 surfaces endowed with a metric fibers over $\gm_{{\rm HK}}$ and the fibre
of $\gm_{{\rm K3}}\to\gm_{{\rm HK}}$ is isomorphic to $S^2$.
Any point in $\gm_{{\rm K3}}$ gives rise to an $N=(2,2)$ SCFT
and the induced inclusion $\gm_{{\rm K3}}\subset \gm_{(2,2)}$ is compatible with the two projections.
Mirror symmetry is realized as a certain discrete group action on $\gm_{(2,2)}$ or $\gm_{(4,4)}$.

Due to the fact that $\gm_{{\rm K3}}\to\gm_{{\rm HK}}$ is only an $S^2$-fibration and not an
$S^2\times S^2$-fibration
as is $\gm_{(2,2)}\to\gm_{(4,4)}$, one soon
realizes that points in $\gm_{{\rm K3}}$ might be mirror symmetric to points that are no longer
in $\gm_{{\rm K3}}$. We will explain that
Hitchin's generalized Calabi--Yau structures allow to give a geometric meaning also to
those points.

From a slightly different point of view, one could think of generalized Calabi--Yau structures
as geometric realizations of points in the extended period domain which is obtained by 
passing  from the period domain $Q\subset \IP(H^2(M,\IC))$, an open subset of a smooth quadric,
to the analogous object $\widetilde Q\subset \IP(H^*(M,\IC))$. The latter
is defined in terms of the Mukai pairing on $H^*(M,\IZ)$.
Recall that due to results of Siu, Todorov, and others,
the period domain $Q$ is essentially the moduli space of marked K3 surfaces. The larger moduli space corresponding
to $\widetilde Q$ contains the $B$-field shifts of those as a hyperplane section. Its complement
is the open subset of $B$-field shifts of symplectic structures on a K3 surface.
Thus, complex structures and symplectic structures are parametrized by the same moduli space and the
discrete group ${\rm O}(H^*(M,\IZ))$ acting on $\widetilde Q$ frequently interchanges these  two.
 
In particular, we will prove the following result:

\begin{theorem}
The period map $\kp_{\rm gen}:\gn_{\rm gen}\to \widetilde Q$
from the moduli space  of generalized Calabi--Yau structures $\gn_{\rm gen}$ on a K3 surface $M$
to the extended period domain $\widetilde Q\subset\IP(H^*(M,\IC))$ is surjective. Moreover, $\gn_{\rm gen}$ admits a 
natural symplectic structure $\Omega$ with respect to which the moduli space of symplectic structures
${\rm Sympl}(M)\subset \gn_{\rm gen}$ is Lagrangian.
\end{theorem}

It might be worth pointing out that the B-field, from a 
mathematical point of view a slightly mysterious object, is indispensable when we want to  view complex and symplectic 
structures as special instances of a more general notion.
I certainly hope and expect that this unified treatment of symplectic and complex structures
on K3 surfaces leads to a better understanding of both. 

\medskip

Here is the plan of the paper. In the first section we recall the notion of generalized Calabi--Yau structures,
which is due to Hitchin, and discuss the two main examples (and their $B$-field transforms) alluded to above.
In Section \ref{GTSect}, after introducing the notion of generalized Calabi--Yau structures of (hyper)k\"ahler type,
we prove a Global Torelli theorem for generalized Calabi--Yau structures on K3 surfaces. We also discuss generalizations
of the existence theorems of Siu and Yau.
Moduli spaces of generalized Calabi--Yau structures are treated in Section \ref{MS}. We define various period
maps and show how they can be used to relate the moduli space of generalized K3 surfaces to the moduli
space of $N=(2,2)$ SCFT.  In Section \ref{Orl}
we argue that these new moduli spaces are well suited to interpret
Orlov's criterion on the equivalence of derived categories
of algebraic K3 surfaces. In order to treat the twisted, still conjectural version of it,
we introduce the Picard group and the transcendental lattice of a generalized Calabi--Yau
structure.
In the last section a natural symplectic (hermitian) structure on the moduli space of generalized Calabi--Yau structures
is defined. It turns out that the part of the moduli space that parametrizes generalized Calabi--Yau structures
of the form $\exp(i\omega)$, with $\omega$ a symplectic form, is Lagrangian.

\medskip

{\bf Acknowledgement:} I am most grateful to Nigel Hitchin; many of the results in this article are directly inspired
by his paper \cite{Hitchin}. I would also like to thank Marco Gualtieri, Klaus Mohnke, 
and Ivan Smith for interesting discussions and helpful comments. 
%%%%%%%%%%%%%%%%%%%%%%%%%%%%%%%%%%%%%%%%%%%%%%%%%%%%%%%%%%%%%%%%%%%%%%%%%%%%%%%%%%%%%%%%%%%%%%%
\section{Hitchin's generalized Calabi--Yau structures}

Throughout this paper we will assume that $M$ is the differentiable manifold underlying
a K3 surface. E.g.\ we could think of $M$ as the differentiable fourfold
defined by $x_0^4+x_1^4+x_2^4+x_3^4=0$ in $\IP^3_\IC$. (Due to a result of Kodaira
one knows that any K3 surface is diffeomorphic to $M$.) We will also
fix the natural orientation induced by the complex structure.
This will enable us to speak about positivity and negativity of four-forms on $M$.

We take the liberty to change some of Hitchin's original  conventions in order to make
the theory compatible with the standard theory of K3 surfaces.

For two even complex forms $\varphi,\psi\in\ka_\IC^{2*}(M)$ one defines 
$$\langle\varphi,\psi\rangle:=-\varphi_0\wedge\psi_4+\varphi_2\wedge\psi_2-\varphi_4\wedge\psi_0\in\ka^4_\IC(M),$$
where $\varphi_i$ and $\psi_i$ denote the parts of degree $i$ of $\varphi$ and $\psi$, respectively. This
is the Mukai pairing on the level of forms.

\begin{definition} A \emph{generalized Calabi--Yau structure} on the four-dimensional
mani\-fold $M$ is a closed even form $\varphi\in\ka^{2*}_\IC(M)$ such that
$$\langle\varphi,\varphi\rangle=0~~{\rm and}~~\langle\varphi,\overline\varphi\rangle>0.$$
\end{definition}

Note that such a $\varphi$ is not necessarily homogeneous and that its degree zero term is constant.

\begin{remark}
Hitchin defines also odd generalized Calabi--Yau structures, but they
are of no importance for our purposes,
as in dimension four they only exist on manifolds with non-trivial first cohomology.
\end{remark}
The notion of generalized Calabi--Yau structures embraces symplectic and complex structures: 

\begin{example}\label{basicexs}
i) Every symplectic structure $\omega$ on $M$ induces a generalized Calabi--Yau structure
$\varphi=\exp(i\omega)=1+i\omega-(1/2)\cdot\omega^2$. 
In order to see that any symplectic structure on $M$ defines the same orientation, i.e.\
that $\omega^2>0$, one can use Seiberg-Witten theory. For our purpose we might as well
just restrict to those.
\silent{Kotschick gives the argument. Roughly, Taubes shows that the Seiberg-Witten invariants
with respect to the reversed orientation are not trivial. On the other hand, the existing
smooth embedded two-spheres get positive square and a result of Fintushel-Stern says that
this implies the triviality of the Seiberg-Witten classes.}

ii) Let  $X$ be a K3 surface. Thus, $X$ is given by a complex structure $I$ on $M$.
The holomorphic two-form $\sigma$, which is unique up to scaling, defines a generalized
Calabi--Yau structure $\varphi=\sigma$.
\end{example}

These two examples are very different from each other  due to the fact that in i)
the constant term is non-trivial, and after scaling we might even assume
that $\varphi_0=1$, whereas the second example $\varphi=\sigma$ has trivial constant term.
In most of the arguments that will follow, one has to distinguish between  these two cases.

\bigskip

If $B$ is a two-form, then $\exp(B)$ acts on $\ka_\IC^*(M)$ by exterior product, i.e.\ 
$$\exp(B)\cdot\varphi=(1+B+(1/2)\cdot B\wedge B)\wedge\varphi.$$
It is easy to see that multiplication with $\exp(B)$ is orthogonal with respect to the
pairing $\langle~~,~~\rangle$, i.e.\
$$\langle\exp(B)\cdot\varphi,\exp(B)\cdot\varphi'\rangle=\langle\varphi,\varphi'\rangle\in\ka^4_\IC(M)$$
for all forms $\varphi,\varphi'$. This immediately yields the following observation due to Hitchin.
\begin{lemma} For any generalized Calabi--Yau structure $\varphi$ and any real closed two-form
$B$, the form $\exp(B)\cdot \varphi$ is again a generalized Calabi--Yau structure.\qed
\end{lemma}

\bigskip

The generalized Calabi--Yau structure $\exp(B)\cdot\varphi$ is called the
\emph{$B$-field transform of $\varphi$}.
Note that $\exp(B)\cdot\exp(i\omega)=\exp(B+i\omega)$.

The following proposition shows that any generalized Calabi--Yau structure is actually a $B$-field transform
of one of the two fundamental examples \ref{basicexs}.

\begin{proposition}\label{Hitcor} (Hitchin) Let $\varphi$ be a generalized Calabi--Yau structure.

{\rm i)} If  $\varphi_0\ne0$, then  $\varphi=\varphi_0\cdot\exp(B+i\omega)$, with $\omega$  a symplectic form and $B$
a closed real two-form.

{\rm ii)} If $\varphi_0=0$, then $\varphi=\exp(B)\cdot \sigma=\sigma+\sigma\wedge B$,
where $\sigma$ is a holomorphic two-form  with respect to some complex structure on $M$ and $B$ is a closed real two-form.
\end{proposition}

\prf i) More explicitly one finds in this case
$$\varphi_0^{-1}\cdot\varphi=\exp\left({\rm Re}(\varphi_0^{-1}\cdot\varphi_2)\right)
\cdot\exp\left(i\cdot{\rm Im}(\varphi_0^{-1}\cdot\varphi_2)\right),$$
Using $\langle\varphi,\overline\varphi\rangle>0$ we obtain
$\varphi_0\overline\varphi_0\left((\overline\varphi_2/\overline\varphi_0)-(\varphi_2/\varphi_0)\right)^2<0$.
Hence, ${\rm Im}(\varphi_0^{-1}\cdot\varphi_2)$ is symplectic. The claimed equality is checked easily. 

ii) Let $\varphi$ be a generalized Calabi--Yau structure with $\varphi_0=0$. In this case $\varphi_2\wedge\varphi_2=0$ and
$\varphi_2\wedge\overline\varphi_2>0$. Due to an observation of Andreotti, there exists a unique complex
structure on $M$ such that $\sigma:=\varphi_2$ is a holomorphic two-form. By definition, the bundle of $(1,0)$-forms
is the kernel of $\varphi_2:\ka_\IC^1\to\ka_\IC^3$. The integrability of the induced almost complex structure is equivalent
to $d\varphi_2=0$. 

Let us first assume that $\varphi_4$ is exact. Any exact four-form can be written
as $\bar\partial\gamma=d\gamma$ for a $(2,1)$-form $\gamma$ (one way to see this is to use
Hodge-decomposition for $\bar\partial$ and the fact that a four-form is exact if and only
if its $d$-harmonic part is trivial if and only if its $\bar\partial$-harmonic
part is trivial). Since $\sigma$ is non-degenerate, there exists a $(0,1)$-form $\delta$ such that
$\sigma\wedge\delta=\gamma$. Clearly, for degree reasons one also knows $\sigma\wedge\bar\delta=0$.
Then with  $B:=d(\delta+\bar\delta)$ one has $B\wedge\sigma=\varphi_4$.

In general, $\varphi_4$ can be written as $\varphi_4=(\varphi_4-\lambda\sigma\overline\sigma)+\lambda\sigma\overline\sigma$ with
$\lambda\in\IC$ such that $\varphi_4-\lambda\sigma\overline\sigma$ is exact.
Then choose a closed form $B'$ with
$B'\wedge\sigma=\varphi_4-\lambda\sigma\overline\sigma$ as before and set $B=B'+\lambda\overline\sigma+\bar\lambda\sigma$.
\silent{Here is another argument which uses the existence of a hyperk\"ahler structure. If $\varphi_4$ is not exact we use the fact that any K3 surface admits a hyperk\"ahler structure
(see the discussion in the later sections). In particular, we can use the harmonic decomposition
of $\varphi_4$ with respect to such a hyperk\"ahler form. Thus, $\varphi_4={\mathcal H}(\varphi_4)+\varphi_4'$
with ${\mathcal H}(\varphi_4)=\lambda(\sigma\wedge\bar\sigma)$ with $\lambda\in\IC$. In particular,
${\mathcal H}(\varphi_4)=\sigma\wedge(\bar\lambda\sigma+\lambda\bar\sigma)$. As was explained before,
the exact four-form $\varphi_4'$ can be written as $\sigma\wedge B'$ for some
real exact two-form $B$. Altogether this shows that $\varphi=\exp(B)\cdot\sigma$, where
$B=\bar\lambda\sigma+\lambda\bar\sigma+B'$ is a real closed two-form.}
\qed
\bigskip

For the notion of isomorphic generalized Calabi--Yau structures we shall 
consider the group ${\rm Diff}_*(M)$ of all diffeomorphisms $f$ of $M$ such that
the induced action $f^*:H^2(M,\IR)\to H^2(M,\IR)$ is trivial. It seems  unknown whether ${\rm Diff}_*(M)$
coincides with the identity component ${\rm Diff}_{\rm o}(M)$ of the diffeomorphism group.
Allowing only ${\rm Diff}_*(M)$ and not the full diffeomorphism group ${\rm Diff}(M)$ might look
not very natural, but for the moduli space considerations it is useful to divide 
by ${\rm Diff}(M)/{\rm Diff}_*(M)={\rm O}^+(M)$ only later (see \cite{Borcea,Donalson2,Huy}).
 
\begin{definition}\label{iso}
Two generalized Calabi--Yau structures $\varphi$ and $\varphi'$ are called \emph{isomorphic} if and only if there
exists an exact $B$-field $B$ and a diffeomorphism $f\in{\rm Diff}_*(M)$ such that
$\varphi=\exp(B)\cdot f^*\varphi'$.
\end{definition}

Clearly, if $\varphi$ and $\varphi'$ are isomorphic generalized
Calabi--Yau structures then $\varphi_0=0$ if and only if $\varphi'_0=0$.
For the two principal examples, i.e.\ $\varphi$ of the form $\sigma$ or $\exp(i\omega)$,
this reduces to the (only slightly modified) standard definition of isomorphisms. Thus,
$\sigma$ and $\sigma'$ are isomorphic as generalized Calabi--Yau structures if and only if there exists
a diffeomorphism $f\in{\rm Diff}_*(M)$ such that $f^*\sigma'=\sigma$. Similarly,
$\exp(i\omega)$ and $\exp(i\omega')$ are isomorphic if and only if there exists
a diffeomorphism $f\in{\rm Diff}_*(M)$ with $\omega=f^*\omega'$.

\begin{remark}
If one wants to build an analogy between symplectic structures $\omega$ and holomorphic forms
$\sigma$, one soon realizes that the natural isotropy groups
${\rm Sympl}(M,\omega)$ and $\{f\in{\rm Diff}(M)~|~f^*\sigma=\sigma\}$ are quite
different in nature. The group of symplectomorphisms is always infinite-dimensional, whereas
$\{f\in{\rm Diff}(M)~|~f^*\sigma=\sigma\}\subset{\rm Aut}(M,I)$ is discrete.
Only when both structures, $\omega$ or rather $\exp(i\omega)$ and $\sigma$,
are considered as generalized Calabi--Yau structures, the analogy emerges:
Namely,
$${\rm Aut}(\varphi=\exp(i\omega))=\{(f,B)~|~\exp(B)\cdot f^*\exp(i\omega)=\exp(i\omega)\}={\rm Sympl}(M,\omega)$$
and
$${\rm Aut}(\varphi=\sigma)=\{(f,B)~|~\exp(B)\cdot f^*\sigma=\sigma\},$$ the  isotropy groups of the generalized
Calabi--Yau structures $\exp(i\omega)$ respectively $\varphi=\sigma$, are both infinite-dimensional.
In fact we see that ${\rm Aut}(\varphi=\sigma)$ is the set of all exact B-fields of type 
$(1,1)$, which can be identified with the space of all functions modulo scalars due 
to the $\partial\bar\partial$-lemma.
\end{remark}

In what follows, we will use the following notation:
\begin{definition} Let $\varphi$ be a generalized Calabi--Yau structure. Then $P_\varphi\subset
\ka^*(M)$ denotes the real vector space spanned by the real and imaginary part of $\varphi$.
Analogously, $P_{[\varphi]}\subset H^*(M,\IR)$ is the plane generated by the real and imaginary
parts of the associated cohomology class.
\end{definition}

Thus, $P_\varphi$ with respect to $\langle~~,~~\rangle$ is positive at every point and
$P_{[\varphi]}\subset H^*(M,\IR)$ is a positive plane with respect to the Mukai pairing.
Moreover, $P_\varphi$ comes along with a natural (pointwise) orientation.
Conversely, the orien\-ted \- plane $P_\varphi\subset\ka^{2*}(M)$ determines $\varphi$ uniquely up to non-trivial complex scalars
(use $d\varphi=0$). Also note that $P_{\exp(B)\varphi}=\exp(B)\cdot P_\varphi$.

Recall that in general there is a natural isomorphism
between the (open subset of a) quadric $Q_V:=\{x~|~x^2=0, x\cdot \bar x>0\}\subset\IP(V_\IC)$
and the Grassmannian of  oriented positive planes ${\rm Gr}_2^{\rm po}(V)$, where $V$ is a real vector space
endowed with a non-degenerate quadratic form.

%%%%%%%%%%%%%%%%%%%%%%%%%%%%%%%%%%%%%%%%%%%%%%%%%%%%%%%%%%%%%%%%%%%%%%%%%%%%%%%%%%%%%%%%%%%%%%%%%
\section{Global Torelli for generalized CY structures of HK type}\label{GTSect}

Recall that any complex structure on $M$ defines a K3 surface (cf.\ \cite{FM}) and therefore
is K\"ahler, due to Siu's results \cite{Siu}. Using the existence of Ricci-flat K\"ahler structures
proved by  Yau in \cite{Yau}, this also shows that any complex structure on $M$ admits
a hyperk\"ahler structure. In fact, any K\"ahler class is represented by a unique hyperk\"ahler
form.

In this section we shall discuss the analogous notions for generalized Calabi--Yau structures.
First note that a symplectic two form $\omega$ is of type $(1,1)$
with respect to a complex structure $I$ if and only if $\sigma\wedge\omega=0$, where $\sigma$ is the holomorphic
two-form on $(M,I)$. In this case, $\omega$ or $-\omega$ is a K\"ahler form.
(As before, we use $\omega^2>0$.)
Thus, if $\omega$ is a symplectic form, then one of the two forms $\omega$ or $-\omega$
is a K\"ahler form with respect to $\sigma$ if and only if $P_\sigma$ and $P_{\exp(i\omega)}$ are pointwise
orthogonal. This can be generalized as follows:

\begin{definition}
Let $\varphi$ be a generalized Calabi--Yau structure on $M$. We say that $\varphi$ is \emph{K\"ahler}
(or of \emph{K\"ahler type})
if there exists another gene\-ralized Calabi--Yau structure $\varphi'$ orthogonal to $\varphi$, i.e.\
such that $P_\varphi$ and $P_{\varphi'}$ are pointwise orthogonal.
In this case, $\varphi'$ is called a \emph{K\"ahler structure} for $\varphi$.
\end{definition}

Note that the orthogonality of two planes $P_\varphi$ and $P_{\varphi'}$ is in general a stronger condition
than just $\langle\varphi,\varphi'\rangle\equiv0$.

%One certainly expects that any generalized Calabi--Yau structure is of K\"ahler type. This
%would be the analogue of Siu's result. This can be shown for generalized Calabi--Yau structures
%$\varphi$ with $\varphi_0=0$ 
Does Siu's existence result of K\"ahler structures on K3 surfaces extend to generalized
Calabi--Yau structures? An affirmative answer can be given for generalized Calabi--Yau structures $\varphi$
with $\varphi_0=0$ (cf.\ Lemma \ref{hk0}).

\begin{example}\label{exk}
i) Let $\varphi=\sigma$. If $\varphi'$ is a K\"ahler structure for $\varphi$, then $\varphi_0'\ne0$,
as $\bigwedge^2$ has only three positive eigenvalues at every point. Thus, we may assume $\varphi'=\exp(B+i\omega)$.
The orthogonality of $P_\varphi$ and $P_{\varphi'}$ is equivalent to $\sigma\wedge B=\sigma\wedge\omega=0$.
Thus, $\varphi'$ is a K\"ahler structure for $\varphi$ if and only if $\varphi'=\exp(B+i\omega)$ (up to scalar factors)
with $B$ a closed real $(1,1)$-form and $\pm\omega$ a K\"ahler form
(both with respect to the complex structure defined by $\sigma$).

ii) Let $\varphi=\exp(i\omega)$, where $\omega$ is a symplectic form, and let $\varphi'$ be a K\"ahler structure
for $\varphi$. There are two possible cases: Either $\varphi'_0=0$, then $\varphi'=\sigma$ and $\pm\omega$ is a K\"ahler form
with respect to the complex structure defined by $\sigma$ or $\varphi'_0\ne0$.
In the latter case $\varphi'=\exp(B'+i\omega')$ (up to scalars). The orthogonality is equivalent to
the four equations $B'\wedge\omega=0$, $B'\wedge\omega'=0$, $\omega\wedge\omega'=0$, and $B'^2=\omega^2+\omega'^2$.
In particular, $\omega$, $\omega'$, and $B'$ are three pairwise pointwise orthogonal symplectic forms.
\silent{
iii) Suppose $\omega$ is a K\"ahler form for $\sigma$. If there exists a diffeomorphism
$f\in{\rm Diff}(M)\setminus{\rm Diff}_*(M)$ such that $f^*\omega$ is still K\"ahler with respect to $\sigma$,
then Moser's result shows that there exists a $g\in {\rm Diff}_*(M)$ with 
$g^*f^*\omega=\omega$. Thus, in such a case one would have $({\rm Diff}(M)\setminus{\rm Diff}_*(M))\cap 
{\rm Sympl}(M,\omega)\ne\emptyset$.}
\end{example}

%{\footnotesize\bf I would like a treatment of K\"ahler structures on the level of the induced generalized complex structures.
%In principle this should be possible, but then I better fix the proof in the previous section.}

Recall that a  K\"ahler form $\omega$ on a K3 surface is a hyperk\"ahler form if $\omega\wedge\omega$ is a scalar
multiple of the canonical volume form $\sigma\wedge\bar\sigma$. Scaling $\sigma$, which does not change the
complex structure, makes it natural to assume that $2\omega\wedge\omega=\sigma\wedge\bar\sigma$.
 
\begin{definition}
A generalized Calabi--Yau structure $\varphi$ is \emph{hyperk\"ahler} if there exists another generalized
Calabi--Yau structure $\varphi'$ such that $\varphi$ and $\varphi'$ are orthogonal and
$\langle\varphi,\bar\varphi\rangle=\langle\varphi',\bar\varphi'\rangle$.
We say that $\varphi'$ is a \emph{hyperk\"ahler structure} for $\varphi$.
\end{definition}

\begin{example}
i) A hyperk\"ahler structure for $\varphi=\sigma$ is a generalized Calabi--Yau structure $\varphi'$
of the form $\lambda\exp(B+i\omega)$, where $0\ne\lambda\in\IC$, $B$ is a closed real $(1,1)$-form,
and $\pm\omega$ is a hyperk\"ahler form such that $2|\lambda|^2\omega\wedge\omega=\sigma\wedge\bar\sigma$.

ii) A hyperk\"ahler structure for $\varphi=\exp(i\omega)$ is either a holomorphic two-form
$\varphi'=\sigma$ with respect to which $\pm\omega$ is a hyperk\"ahler form or
it is of the form $\varphi'=\exp(B'+i\omega')$ (up to scalar factors which we omit) as in ii) of Example \ref{exk}
with the additional condition $\omega\wedge\omega=\omega'\wedge\omega'$. This shows that
$\sigma:=(1/\sqrt2)B'+i\omega'$ defines a complex structure with respect to which $\pm\omega$ is a hyperk\"ahler form.
\end{example}

\begin{remark} i) Clearly, both definitions are symmetric in $\varphi$ and $\varphi'$, i.e.\
if $\varphi'$ is a (hyper)k\"ahler structure for $\varphi$ then
$\varphi$ is a (hyper)k\"ahler structure for $\varphi'$.

ii) Let $\varphi$ be a generalized Calabi--Yau structure and $\varphi'$ a (hyper)k\"ahler
structure for it. Then $\exp(B)\cdot\varphi'$ is a (hyper)k\"ahler structure for the $B$-field
transform $\exp(B)\cdot\varphi$.

\silent{iii) Note that if $\varphi'$ and $\varphi''$ are isomorphic generalized Calabi--Yau structures
and $\varphi'$ is a (hyper)k\"ahler structure for a generalized Calabi--Yau structure
$\varphi$ then $\varphi''$ does not necessarily define a (hyper)k\"ahler structure for
$\varphi$.  But note that if $\varphi$ and $\varphi'$ are given, then
there exists a generalized Calabi--Yau structure $\tilde\varphi$  isomorphic to $\varphi$ and defining a
(hyper)k\"ahler structure for $\varphi'$ if and only if there exists a generalized Calabi--Yau
structure $\tilde\varphi'$ isomorphic to $\varphi'$ and defining a (hyper)k\"ahler structure for
$\varphi$. Indeed, if $\tilde\varphi=\exp(B)\cdot f^*\varphi$ then set
$\tilde\varphi'=\exp(-f^*B)f^*\varphi'$.}
\end{remark}

Obviously, any generalized Calabi--Yau structure which is hyperk\"ahler is also K\"ahler.
The following lemma thus settles the existence question for both structures 
in the case $\varphi_0=0$.

\begin{lemma}\label{hk0}
Any generalized Calabi--Yau structure $\varphi$ with $\varphi_0=0$ is hyperk\"ahler.
\end{lemma}

\prf As we have seen, a generalized Calabi--Yau structure
$\varphi$ with $\varphi_0=0$ is of the form $\sigma+\sigma\wedge B$, where
$\sigma$ is a holomorphic two-form with respect to a certain complex structure $I$ on $M$.
Using the results of Siu and Yau we find a  hyperk\"ahler form $\omega$ on $(M,I)$.
Thus, $\exp(i\omega)$ defines a hyperk\"ahler structure for $\sigma$.
Using the above remark, we find that $\varphi'=\exp(B+i\omega)$ is a hyperk\"ahler structure
for $\varphi=\exp(B)\cdot\sigma$.\qed

\bigskip

Recall that Yau's existence result says that for any
complex structure defined by a complex two-form $\sigma$ and any
K\"ahler form $\omega$ there exists a hyperk\"ahler form
$\omega'$ cohomologous to $\omega$. Using Moser's result \cite{Moser}, which shows that $\omega$ and
$\omega'$ are related by a diffeomorphism $f\in {\rm Diff}_{\rm o}(M)$, we find that Yau's result is equivalent
to saying that any K\"ahler structure for the generalized Calabi--Yau structure
$\sigma$ is isomorphic to a hyperk\"ahler structure.

\begin{proposition}
Let $\varphi'$ be a K\"ahler structure for a generalized Calabi--Yau structure
$\varphi$. If $\varphi_0\cdot\varphi'_0=0$ then $\varphi'$ is isomorphic to a 
hyperk\"ahler structure for $\varphi$.
\end{proposition}

\prf
Let us first assume $\varphi_0=0$.  Then $\varphi'$ is necessarily (up to scaling) of the form $\exp(B+i\omega)$.
The assertion is invariant under shifting both structures by a $B$-field. Thus, we may assume
that $\varphi=\sigma$. Yau's result immediately proves the existence of a diffeomorphism
$f\in {\rm Diff}_*(M)$ such that $f^*\omega$ is a hyperk\"ahler form for $\sigma$.
Moreover, $f^*B$ and $B$ differ by an exact $B$-field $B'$, i.e.\ $B-f^*B=B'$.
Therefore, $\exp(B')\cdot f^*\exp(B+i\omega)$ is a hyperk\"ahler structure (up to scalar factors)
for the generalized Calabi--Yau structure $\sigma$. 

Next assume $\varphi_0\ne0$.
After rescaling we have $\varphi=\exp(B+i\omega)$,
where $\omega$ is a symplectic form. Clearly, $\varphi$ is (hyper)k\"ahler if and only if $\exp(i\omega)$ is
(hyper)k\"ahler. Thus, we may assume $\varphi=\exp(i\omega)$.
By assumption $\varphi'$ is a K\"ahler structure  for
 $\varphi$. A priori, we have to distinguish the two cases $\varphi'_0=0$ and $\varphi'_0\ne0$, but the
 second case is excluded by assumption.

If $\varphi'_0=0$, then $\varphi'=\sigma+\sigma\wedge B$. The orthogonality of $P_\varphi$ and $P_{\varphi'}$ 
yields $\sigma\wedge B=0$, i.e.\ $\varphi'=\sigma$. Thus, $\pm\omega$ is a K\"ahler structure with respect to
the complex structure defined by $\sigma$. Hence, there exists a (unique) hyperk\"ahler form $\omega'$
cohomologous to $\omega$ which can in fact be written
as $\omega'=f^*\omega$ for some $f\in{\rm Diff}_{\rm o}(M)$ due to the result
of Moser.  But then $\pm\omega$ is a hyperk\"ahler form with respect to the complex structure defined by
$(f^{-1})^*\sigma$ (up to sign).
Hence, $(f^{-1})^*\sigma$ is a hyperk\"ahler structure for $\exp(i\omega)$ which is isomorphic to $\varphi'$
via $f$.
\qed

\begin{remark}
I certainly believe that the hypothesis $\varphi_0\cdot\varphi'_0=0$ is superfluous.
The problem is a situation where $\varphi=\exp(i\omega)$ and $\varphi'=\exp(B'+i\omega')$.
The K\"ahler condition is equivalent to
$B'\wedge\omega=\omega'\wedge\omega=B'\wedge \omega'=0$ and $B'^2=\omega^2+\omega'^2$.
The problem one has to solve in this case seems very similar to the original existence question for
Ricci-flat K\"ahler forms.

In order to see this analogy, we suppose for simplicity that $\int\omega^2=\int\omega'^2$.
Then consider the complex two-form $\sigma:=\omega+i\omega'$  which is clearly
orthogonal to $B'$. Moreover, $B'$ satisfies $B'^2=\sigma\overline\sigma$.
If we can change $B'$ by an exact form, such that the new $B'$ is still orthogonal to $\sigma$ and 
$B'^2=2{\rm Im}(\sigma)^2$, then we are done. Indeed, then $\sigma':=B'+i\sqrt{2}\omega'$ would be
a Calabi--Yau form orthogonal to $\omega$. The latter could be made a hyperk\"ahler form by applying
a diffeomorphism $f\in{\rm Diff}_*(M)$ due to Yau's theorem.
Note that $\sigma$ itself is a priori not a Calabi--Yau form, since, $\sigma^2\ne0$ {\it a
priori} but the "$(1,1)$"-form $B$
satisfies already the condition $B^2=\sigma\overline\sigma$.
\end{remark}

\begin{remark}\label{donaldsonansatz}
It is an open question whether any symplectic form on $M$ is in fact (hyper)k\"ahler
with respect to some complex structure. One expects an affirmative answer to this
and a possible approach has recently been suggested by Donaldson \cite{Donaldson}.
The last proposition extended to the case $\varphi_0\cdot \varphi_0'\ne0$ would 
show in particular that if the generalized Calabi--Yau structure $\exp(i\omega)$
associated to a symplectic form $\omega$ is K\"ahler (as a generalized Calabi--Yau structure),
then $\omega$ is actually a
hyperk\"ahler form with respect to a certain complex structure on $M$.
Thus, together with an analogue of Siu's existence
result, which would claim that any generalized Calabi--Yau structure is of K\"ahler 
type, the more general version of the above proposition would in particular show that any symplectic form
is hyperk\"ahler.
\end{remark}

\begin{remark}
When conjecturing the existence of Ricci-flat metrics, Calabi gave a simple proof of the unicity, i.e.\
any K\"ahler class is represented by at most one hyperk\"ahler form. Equivalently,
if $f\in{\rm Diff}_*(M)$ such that $f^*\sigma=\sigma$, then also
$f^*\omega=\omega$ for any hyperk\"ahler form $\omega$ on the complex K3 surface determined by $\sigma$.
Thus, $f$ is an isometry and hence of finite order. It is known that this in fact
yields $f={\rm id}$. The unicity is no longer true when the role of $\sigma$
and $\omega$ are interchanged,  i.e.\ for a given $\omega$ there may
exist several complex structures realizing the same period and making $\omega$ a hyperk\"ahler form.
Indeed, if  $\omega$ is a hyperk\"ahler form with respect to $\sigma$ and ${\rm id}\ne f\in
{\rm Sympl}(\omega)\cap{\rm Diff}_*(M)$,
then $f^*\sigma\ne\sigma$ by the above argument.
Hence, $\sigma$ and $f^*\sigma$ are two different hyperk\"ahler structures for $\exp(i\omega)$. 

Again, the different behaviour of $\omega$ and $\sigma$ can be explained if both are considered
as generalized Calabi--Yau structures. Indeed, for $\varphi=\sigma$ there exist many different hyperk\"ahler structures
$\exp(B+i\omega)$ in the same cohomology class. In fact, if $\omega$ is hyperk\"ahler for $\sigma$, then
$\exp(B+i\omega)$ is a hyperk\"ahler generalized Calabi--Yau structure for $\varphi=\sigma$ whenever
$B$ is  an exact $(1,1)$-form. 
\end{remark}

The arguments used in the proof of the following result
show in particular that two hyperk\"ahler structures for a given generalized Calabi--Yau structure
are always isomorphic.

\begin{proposition}\label{GT}(Global Torelli theorem)
Let $\varphi$ and $\psi$ be two gene\-ra\-lized Calabi--Yau structures on $M$  and suppose they are both hyperk\"ahler.
If  $P_{[\varphi]}=P_{[\psi]}\subset
H^*(M,\IR)$, then there exists a real exact $B$-field $B$ and a diffeomorphism
$f$ such  that $P_\varphi=\exp(B)\cdot P_{f^*\psi}$, i.e.\ up to rescaling $\varphi=\exp(B)\cdot f^*\psi$.

If $\varphi_0\ne0$, then $f$ can be chosen in ${\rm Diff}_*(M)$, i.e.\ $\varphi$ and $\psi$ are isomorphic
generalized Calabi--Yau structures.
\end{proposition}

\prf First suppose that $\varphi_0=0$. Then also $\psi_0=0$ and after rescaling we may assume
that $[\varphi_4]=[\psi_4]$. We have to find a real exact two-form $B$ and a diffeomorphism $f$
such that $\varphi_2=f^*\psi_2$ and $\varphi_4=B\wedge f^*\psi_2+f^*\psi_4$.
As has been explained before, the assumption that $\varphi$ and $\psi$ are generalized
Calabi--Yau structures implies that $\varphi_2$ and $\psi_2$ are holomorphic two-forms
with respect to uniquely determined complex structures. Invoking the classical Global Torelli theorem
for K3 surfaces we find a diffeomorphism $f$ such that $\varphi_2=f^*\psi_2$.
Thus, we may assume $\varphi_2=\psi_2$ already and try to find an real exact  two-form
$B$ such that $\varphi_4-\psi_4=B\wedge\psi_2$. This follows directly from the argument given in the
proof of Proposition \ref{Hitcor}, because $\varphi_4-\psi_4$ is exact and $\psi_2$ is a non-degenerate holomorphic two-form.
(Note that $f$ can be chosen in ${\rm Diff}_*(M)$ if $\varphi$ and $\psi$ admit hyperk\"ahler structures
$\varphi'$ respectively $\psi'$ with $[\varphi']=[\psi']$.)

If $\varphi_0\ne0$ then also $\psi_0\ne0$ and after rescaling we might
assume $\varphi_0=\psi_0=1$. Thus, we have $\varphi=\exp(B+i\omega)$
and $\psi=\exp(B'+i\omega')$, where $\omega$ and $\omega'$ are symplectic forms
and $B$ and $B'$ are cohomologous real closed two-forms. In particular,
$B$ and $B'$ differ by an exact $B$-field. Thus, we may reduce
to the case $\varphi=\exp(i\omega)$ and $\psi=\exp(i\omega')$ with $[\omega]=[\omega']$.

The assumption ensures that $\omega$ and $\omega'$ are hyperk\"ahler forms
with respect to certain holomorphic two-forms $\sigma$ and $\sigma'$, respectively.

As the moduli space of K3 surfaces (with metric) is connected, there exists
a deformation $(\omega_t,\sigma_t)$ of $(\omega_0,\sigma_0)=(\omega,\sigma)$
such that $([\omega_1],[\sigma_1])=([\omega'],[\sigma'])$.
Moreover, since $[\omega]=[\omega']$, we may assume that $[\omega_t]$ is constant.
Using Moser's result we then find a continuous family of diffeomorphisms $f_t$ such that
$f_t^*\omega_t=\omega$. Applying $f$ also to $\sigma_t$ shows that we can in fact
assume that $(\omega_t,\sigma_t)\equiv(\omega,\sigma_t)$. The upshot of all this is that whenever $\varphi=\exp(i\omega)$
and $\psi=\exp(i\omega')$ are two hyperk\"ahlerian generalized Calabi--Yau structures, then
we can choose $\sigma$ and $\sigma'$ such that $[\sigma]=[\sigma']$. The standard Global Torelli theorem
then yields the existence of a diffeomorphism $f\in{\rm Diff}_*(M)$ such that
$f^*\sigma'=\sigma$.
The latter condition and the unicity of the Ricci-flat K\"ahler form representing a given
K\"ahler class show $f^*\omega'=\omega$.
\qed

\begin{remark}
i) We can slightly improve the above statement. Assume $\varphi$ and $\psi$ are two generalized
Calabi--Yau structures. Suppose that there exists an automorphism $F$ of
the K3 lattice $H^2(M,\IZ)$  such that $F[\psi]=[\varphi]$. (Here, we extend
$F$ by the identity to the full cohomology.) Then there exists
a diffeomorphism $f\in{\rm Diff}(M)$ and an exact two-form $B$ such that
$\varphi=(\exp(B)\circ(\pm{\rm id}_{H^2})\circ f^*)(\psi)$ and 
$f^*=F$.  Indeed, due to a result of Borcea and Matumuto there exists a diffeomorphism $F$ such that
$f^*=\pm F$. Then we may apply the proposition to $\varphi$ and $f^*\psi$. 

ii) A Global Torelli theorem for generalized Calabi--Yau structures $\varphi$ and $\psi$
which are not necessarily hyperk\"ahler would in particular show that for any two cohomologous
symplectic structures $\omega$ and $\omega'$ on $M$ there exists a diffeomorphism $f$ such that
$\omega=f^*\omega'$. This could in turn be used to show that every symplectic structure on $M$ is hyperk\"ahler
(cf.\ Remark \ref{donaldsonansatz}).
\end{remark}

%%%%%%%%%%%%%%%%%%%%%%%%%%%%%%%%%%%%%%%%%%%%%%%%%%%%%%%%%%%%%%%%%%%%%%%%%%%%%%%%%%%%%%%%%%%%%%%%%%
\section{Generalized K3 surfaces and moduli spaces}\label{MS}

Marked K3 surfaces endowed with a K\"ahler structure and a $B$-field form a moduli space
that can be described via the period map. It turns out that the period map injects this moduli
space into the physics moduli space of $N=(2,2)$ SCFT. However, not every $N=(2,2)$ SCFT parametrized
by the latter comes from a classical K3 surface (with a $B$-field);  the geometric moduli space is of real
codimension two.
In fact, the $N=(2,2)$ SCFT moduli space fibers over the $N=(4,4)$ SCFT moduli space with fibre
$S^2\times S^2$. The standard K3 moduli space fibers as well over the  $N=(4,4)$ SCFT moduli space
(which is interpreted as the moduli space
of hyperk\"ahler metrics), but the fibre is only $S^2$, the twistor line.

In this section we shall indicate how Hitchin's generalized Calabi--Yau structures (or rather
generalized K3 surfaces) fit nicely in this picture. We will see that
the $N=(2,2)$ SCFT moduli space can be interpreted as the moduli space of generalized K3 surfaces.
For details of the moduli space construction we refer the reader to \cite{Huy} or the original articles \cite{AM,NW}.
In particular, we will use the notations introduced there. Also note that we will actually not work with
the SCFT moduli spaces, but rather with certain period domains, which have
been shown to contain the corresponding moduli spaces (cf.\ \cite{AM}).

\bigskip

So far, by a K3 surface we meant a compact complex  surface $X$ with trivial
canonical bundle $K_X$ and $b_1(X)=0$. Any K3 surface in this sense is determined by
a complex structure $I$ or by a Calabi--Yau structure $\sigma\IC$ on $M$.
From now on we will reserve the name K3 surface for a complex
surface already endowed with a hyperk\"ahler form. More precisely, we have

\begin{definition} A \emph{K3 structure} on the differentiable manifold $M$
consist of a closed complex two-form $\sigma\in\ka^2_\IC(M)$ with $\sigma\wedge\sigma=0$ and
a symplectic
form $\omega\in\ka^2(M)$ such that {\rm i)} $\omega\wedge\sigma=0$ and
{\rm ii)} $\sigma\wedge\bar\sigma=2\omega\wedge\omega>0$.
\end{definition}

As was explained before, the complex two-form $\sigma$ defines a unique complex structure. The orthogonality
condition $\sigma\wedge\omega=0$ is equivalent to $\omega$ being a $(1,1)$-form with respect to this
complex structure. Eventually, $\sigma\wedge\bar\sigma=2\omega\wedge\omega$ ensures that $\pm\omega$ is a hyperk\"ahler
form. Thus, $M$ endowed with such a K3 structure is just a K3 surface with a chosen hyperk\"ahler structure.
Using the convention of the last section we give the following
 
\begin{definition}\label{defiso} A \emph{generalized K3 structure} on $M$ consists of a pair $(\varphi,\varphi')$
of generalized Calabi--Yau structures $\varphi$ and $\varphi'$ such that $\varphi$ is a hyperk\"ahler
structure for $\varphi'$.

Two generalized K3 structures $(\varphi,\varphi')$ and $(\psi,\psi')$ on $M$ are called
isomorphic if there exists a diffeomorphism $f\in{\rm Diff}_*(M)$ and an exact real two-form
$B\in\ka^2(M)$ such that $(\varphi,\varphi')=\exp(B)\cdot f^*( \psi,\psi')$.
\end{definition}

Clearly, the $B$-field transform $(\exp(B)\cdot\varphi,\exp(B)\cdot\varphi')$
of any  generalized K3 structure $(\varphi,\varphi')$ is again a generalized K3 structure.

\begin{definition} To any generalized K3 structure $(\varphi,\varphi')$ on $M$ we associate
the (pointwise) oriented positive 
four-space $\Pi_{(\varphi,\varphi')}\subset \ka^{2*}(M)$ spanned by $P_\varphi$ and $P_{\varphi'}$.
Analogously, one defines an oriented positive four-space $\Pi_{([\varphi],[\varphi'])}\subset H^*(M,\IR)$
spanned by  $P_{[\varphi]}$ and $P_{[\varphi']}$.
\end{definition}

\begin{remark}
The set $T_\Pi$ of all generalized K3 structures $(\varphi,\varphi')$ with fixed positive four-space $\Pi$
is naturally isomorphic to the Grassmannian of oriented planes ${\rm Gr}^{\rm o}_2(\Pi)=S^2\times S^2=\IP^1\times\IP^1$.
Indeed, $T_\Pi=Q_\Pi$, which is a quadric in $\IP^3=\IP^3(\Pi_\IC)$. We call $T_\Pi$ or $T_{(\varphi,\varphi')}$ the
(generalized) \emph{twistor space} (or, more precisely, the base of it).
\end{remark}

\begin{example} If $(\varphi=\sigma,\varphi'=\exp(i\omega))$ is a classical K3 structure on $M$,
then $\Pi$ is spanned by the oriented base ${\rm Re}(\sigma),{\rm Im}(\sigma),1-(1/2)\cdot\omega^2=1-(1/4)\sigma\bar\sigma,
\omega$. In other words, if we write $\sigma=\omega_J+i\omega_K$, where $J,K=IJ$ are the two other natural
complex structures induced by the hyperk\"ahler form, then
$\Pi=\langle 1-(1/2)\omega^2,\omega_I=\omega,\omega_J,\omega_K\rangle$.
The classical twistor deformations $S^2=\{aI+bJ+cK~|~a^2+b^2+c^2=1\}$ form a $\IP^1$ which is contained
in $T_{(\sigma,\exp(i\omega))}$ as the hyperplane section $Q_\Pi\cap \IP(\langle\omega_I,\omega_J,\omega_K\rangle_\IC)$.
In particular, it is not one of the two components of $\IP^1\times\IP^1$.

Note that the other generalized K3 structures parametrized by
$\IP^1\times\IP^1\setminus S^2$ are not obtained as $B$-field transforms of points in $S^2$.

\end{example}

\begin{proposition}\label{clgen}
Let $(\varphi,\varphi')$ be a generalized K3 structure. Then there exists a classical K3 structure
$(\sigma,\exp(i\omega))$  and a closed $B$-field $B$ with
 $$\Pi_{(\varphi,\varphi')}=\exp(B)\cdot\Pi_{(\sigma,\exp(i\omega))}.$$
\end{proposition}
\prf
Since $\Pi_{([\varphi],[\varphi'])}\subset H^*(M,\IR)$ is a positive four-space and
$H^0\oplus H^4$ is only two-dimensional, there exists a positive plane
$H\subset\Pi_{([\varphi],[\varphi'])}\cap H^2(M,\IR)$.
\silent{(Indeed, if $\alpha,\beta,\gamma,\delta$ is an orthogonal basis of $\Pi_{([\varphi],[\varphi'])}$, then we may 
change the orthogonal basis of $\langle\alpha,\beta\rangle$ and $\langle\gamma,\delta\rangle$ such that
$\alpha_0=0$ and $\gamma_0=0$. After this we can modify the orthogonal basis
of $\langle\beta,\delta\rangle$ such that also $\beta_0=0$. In $\langle\alpha,\beta,\gamma\rangle$ we then
find two orthogonal vectors $\alpha,\beta$ with
$\alpha_0=\alpha_4=\beta_0=\beta_4=0$.)}

Using the isomorphism $\Pi_{(\varphi,\varphi')}\cong\Pi_{([\varphi],[\varphi'])}$,
we may choose a generalized Calabi--Yau structure $\psi\in\ka^{2*}_\IC(M)$ 
such that $\langle{\rm Re}(\psi),{\rm Im}(\psi)\rangle
\subset\Pi_{(\varphi,\varphi')}$ corresponds to $H$. 
Hence, $\psi$ is of the form $\exp(B)\cdot\sigma$ with $B$ exact (cf.\ Prop. \ref{Hitcor}).

The orthogonal plane $H^\perp\subset\Pi_{(\varphi,\varphi')}$ is spanned by real and imaginary part of a 
form of type $\exp(B'+i\omega)$ for some closed $B$-field $B'$ and a symplectic structure $\omega$.
(We use here that the Mukai pairing on $H^2(M,\IR)$ has only three positive eigenvalues.) Then, $(\psi_2,\exp(i\omega))$
is automatically a classical K3 structure and, moreover, $\psi_2$ and $B'-B$ are orthogonal.
Thus, $\Pi_{(\varphi,\varphi')}=\Pi_{(\exp(B)\cdot\psi_2,\exp(B'+i\omega))}=\exp(B')\cdot\Pi_{(\psi_2,\exp(i\omega))}$.
\qed

\begin{remark}
The relation between generalized K\"ahler structures and so-called bi-hermitian
structures is explained in detail in \cite[Sect.6]{Gualtieri}. Gualtieri's notion
of generalized K\"ahler structures is formulated in terms of generalized complex structures
(and not generalized Calabi--Yau structures as in this text) and is in fact slightly stronger
than ours.
\end{remark}

\begin{definition} Let $$\gm_{\rm genK3}=
\raisebox{.0ex}{$\{(\varphi,\psi)\}$}_/\raisebox{-1.30ex}{$\cong$}$$ be the
\emph{moduli space of all generalized K3 structures} modulo isomorphism as defined in \ref{defiso}.
\end{definition}

As a consequence of the above proposition \ref{clgen} we obtain
\begin{corollary}
There exists a natural $S^2\times S^2$-fibration 
$$\xymatrix{\gm_{\rm genK3}\ar@{->>}[r]& \gm_{{\rm HK}}:=\left(\raisebox{.6ex}{${\rm Met}^{\rm HK}(M)$}/
\raisebox{-.6ex}{${\rm Diff}_*(M)$}\right)\times H^2(M,\IR)}$$
onto the space of all $B$-field shifts of hyperk\"ahler metrics on $M$.\qed
\end{corollary}

%If $(\varphi,\varphi')$ and $(\psi,\psi')$ are isomorphic via $(B,f)$ then 
%$([\psi],[\psi'])=([\varphi],[\varphi'])$. The converse also holds true:

\begin{definition} The \emph{period} of a generalized K3 structure $(\varphi,\psi)$ is the orthogonal
pair of positive oriented planes $(P_{[\varphi]},P_{[\psi]})\in{\rm Gr}_{2,2}^{\rm po}(H^*(M,\IR))$.
The \emph{period map}
$$\xymatrix{\kp_{\rm genK3}:\gm_{\rm genK3}\ar[r]&{\rm Gr}_{2,2}^{\rm po}(H^*(M,\IR))}$$
is the map that associates to a generalized K3 structure
$(\varphi,\psi)$ its period $(P_{[\varphi]},P_{[\psi]})$.
\end{definition}

%Sometimes, also the positive oriented four-space $\Pi_{([\varphi],[\psi])}=P_{[\varphi]}\oplus P_{[\psi]}
%\in{\rm Gr}_{4}^{\rm po}(H^*(M,\IR))$
%is called the period  of the generalized K3 structure, but recall that the
%map $${\rm Gr}_{2,2}^{\rm po}(H^*(M,\IR))\to{\rm Gr}_{4}^{\rm po}(H^*(M,\IR))$$ is an
%$S^2\times S^2$-fibration.

Similarly, one defines the period map
$$\gm_{{\rm HK}}\rpfeil{9}{\kp_{\rm HK}}{\rm Gr}_{4}^{\rm po}(H^*(M,\IR)),$$
which maps the $B$-field shift of a hyperk\"ahler metric given by a generalized K3 structure
$(\varphi,\psi)$ to $\Pi_{([\varphi],[\psi])}$. We obtain a commutative diagram, where both vertical arrows
are $S^2\times S^2$ fibrations:

$$
%\begin{array}{ccc}
\xymatrix{\gm_{\rm genK3}\ar[d]\ar[rr]^-{\kp_{\rm genK3}}&&{\rm Gr}_{2,2}^{\rm po}(H^*(M,\IR))\ar[d]\\
\gm_{{\rm HK}}\ar[rr]^-{\kp_{\rm HK}}&&{\rm Gr}_{4}^{\rm po}(H^*(M\IR))}
%\end{array}
$$

Classically, one considers the space
$\gm_{\rm K3}$ of  pairs $(I,\omega)$, where $I$ is a complex structure and $\omega$ is a hyperk\"ahler
form modulo the group ${\rm Diff}_*(M)$. This is equivalent to giving the holomorphic
two-form $\sigma$ (up to scaling)
and the K\"ahler form $\omega$.
Adding the $B$-field one obtains the moduli space $\gm_{\rm K3}\times H^2(X,\IR)$.
The fibre of the natural map
$\gm_{\rm K3}\times H^2(X,\IR)\to \gm_{{\rm HK}}$, which
associates to $((\sigma,\omega),B)$ the underlying hyperk\"ahler metric $(g,B)$ (everything shifted by $B$),
is a natural $S^2$-bundle, where
every fibre parametrizes all complex structures associated with one hyperk\"ahler metric.

Using the exponential map we obtain a canonical injection
$$\gm_{\rm K3}\times H^2(X,\IR)\hookrightarrow\gm_{\rm genK3},~~(\sigma,\omega,B)\mapsto (\exp(B)\cdot\sigma,\exp(B+i\omega)).$$
Both sides are fibred over $\gm_{{\rm HK}}$ with
fibre $S^2$ and $S^2\times S^2$, respectively.

\begin{proposition}
The period map
$$\xymatrix{\kp_{\rm genK3}:\gm_{\rm genK3}\ar[r]&{\rm Gr}_{2,2}^{\rm po}(H^*(M,\IR))}$$
is an immersion with dense image.
\end{proposition}

\prf
This is essentially a consequence of the known results for classical K3 surfaces.
One knows that ${\kp_{\rm HK}}:\gm_{{\rm HK}}\to{\rm Gr}_{4}^{\rm po}$ is a dense immersion.
Thus, the assertion follows from the fact that both maps
${\kp_{\rm genK3}}:\gm_{\rm genK3}\to\gm_{{\rm HK}}$ and
${\rm Gr}^{\rm po}_{2,2}\to{\rm Gr}^{\rm po}_4$ are $S^2\times S^2$-fibrations, whose fibres over points
of the image of $\kp_{\rm HK}$ are naturally identified via $\kp_{\rm genK3}$.\qed

\bigskip

The non-surjectivity is a phenomenon already encountered
on the level of hyperk\"ahler metrics. The period map 
$\gm_{{\rm HK}}\rpfeil{5}{}{\rm Gr}_{4}^{\rm po}(H^*(M,\IR))$ is not surjective, but
its image is dense and  points in the complement can be interpreted as
degenerate hyperk\"ahler metrics. (In fact, only positive four-spaces which are not orthogonal
to any $(-2)$-class should be interpreted in this way and should be considered as points defining
a $N=(4,4)$ SCFT.)

\bigskip
So far our discussion of the relation between generalized K3 surfaces and $N=(2,2)$ SCFT
was purely on the level of moduli spaces.
Leaving aside the conformal aspects of the theory, which in any case are not
mathematically established for K3 surfaces, and considering only the supersymmetric part, one can
in fact show that generalized K3 surfaces are indeed related to $N=(2,2)$ supersymmetry.

Let us  briefly recall the general setting for supersymmetry in (hyper)k\"ahler geome\-try. If $X$ is
a compact K\"ahler manifold endowed with a K\"ahler form $\omega$, then the Lefschetz operator
$L$ and its dual $\Lambda$ generate
an $\ssl(2,\IC)$ subalgebra of the algebra of endomorphisms of $\ka^*_\IC(X)$. Together with the
differential operators $d$ and $d^*$ they generate a Lie algebra which in physics jargon
is called the $N=(2,2)$ supersymmetry algebra associated to $(X,\omega)$. Note that this Lie algebra
is finite dimensional and that as a vector space it is generated by $\partial,\bar\partial,
\partial^*,\bar\partial^*,\Delta, L,\Lambda, H$. (The Lie algebra `closes'.) The K\"ahler condition
is crucial at this point. 

If $X$ is a K3 surface and $\omega$ is Ricci-flat one can find an even
bigger Lie algebra. Firstly, as Verbitsky \cite{Verb} has shown more generally for hyperk\"ahler mani\-folds, the
Lefschetz ope\-rators $L_I,L_J,L_K$ and their dual $\Lambda_I,\Lambda_J,\Lambda_K$ asso\-ciated to the
three complex structures $I,J,K=IJ$ generate a Lie subalgebra isomorphic to $\sso(5,\IC)$.
The $N=(4,4)$ supersymmetry algebra associated with $(X,\omega)$  is the (finite dimensional!)
Lie algebra generated by this $\sso(5,\IC)$ and $d$ and $d^*$.
As before, it also contains the differential operators $\bar\partial_I, \bar\partial_J,\bar\partial_K$,
etc. Note that the $N=(4,4)$ supersymmetry algebra is naturally induced by the underlying
hyperk\"ahler metric.

For K3 surfaces, the description of the Lie algebra generated by the Lefschetz operators and their dual
is rather straightforward. Consider the five-dimensional vector space
$V$ spanned by $1\in\ka^0(X),\omega_I,\omega_J, \omega_K\in\ka^2(X),\omega_I^2\in\ka^4(X)$
endowed with the Mukai pairing. Clearly, all six
Lefschetz operators preserve this space and it can be checked that they are all in $\sso(V)$.
Note that singling out a specific complex structure $\lambda=aI+bJ+cK$ naturally yields an inlcusion
$\ssl(2,\IC)\subset\sso(5,\IC)$. In particular, to any such $\lambda$ the $N=(2,2)$ supersymmetry
algebra associated with $\lambda$ and $\omega_\lambda$ is naturally contained in the $N=(4,4)$ 
supersymmetry algebra associated with the hyperk\"ahler metric.

More generally, each point in the fibre of $\gm_{(2,2)}\to\gm_{(4,4)}$, the map that sends
a $N=(2,2)$ SCFT to the induced $N=(4,4)$ SCFT, should in parti\-cular single out a $N=(2,2)$ supersymmetry algebra
within the $N=(4,4)$ supersymmetry algebra determined by a point in $\gm_{(4,4)}$. The following proposition
shows that also the fibre of $\gm_{\rm genK3}\to\gm_{\rm HK}$ naturally parametrizes $N=(2,2)$ supersymmetry
algebras within the $N=(4,4)$ supersymmetry algebra given
by a hyperk\"ahler metric $g\in\gm_{\rm HK}$.

Let $g$ be a hyperk\"ahler metric. Denote by $\Pi_g\subset\ka^{2*}(M)$ and
$\gggg_g\cong\sso(5,\IC)\subset{\rm End}(\ka^{2*}_\IC(M))$ the positive four-space 
respectively the Lie subalgebra 
naturally associated with $g$. Recall that giving a generalized Calabi--Yau structure
$\varphi$ with $P_\varphi\subset\Pi_g$ is equivalent to giving a generalized
K3 surface $(\varphi,\varphi')$ realizing $\Pi_g$.

\begin{proposition}
Any generalzed Calabi--Yau structure $\varphi$ on $M$ with $P_\varphi\subset\Pi_g$
determines a Lie subalgebra of $\gggg_g\cong\sso(5,\IC)$ naturally isomorphic to $\ssl(2,\IC)$.
\end{proposition}

\prf We know that $\Pi_g=\langle1-\omega_I^2/2,\omega_I,\omega_J,\omega_K\rangle_\IR$ and 
$\gggg_g=\sso(V_\IC)$ where $V=\langle1,\omega_I,\omega_J,\omega_K,\omega_I^2\rangle_\IR$. The standard 
isomorphism $\ssl(2,\IC)\cong\sso(3,\IC)$ can be interpreted in our context as $\ssl( (P_\varphi)_\IC)
\cong\sso((P_\varphi\oplus {\Pi_g}^\perp)_\IC)$.  Since $V=P_\varphi\oplus \Pi_g^\perp\oplus P_\varphi^\perp$,
where $P_\varphi^\perp$ is the orthogonal complement of $P_\varphi\subset \Pi_g$,
we obtain a natural inclusion $\ssl(2,\IC)\cong\ssl({P_\varphi}_\IC)\subset \sso(V_\IC)\cong\sso(5,\IC)$.

We leave it to the reader to verify that for $\varphi=\exp(i\omega_\lambda)$ the inclusion 
$\ssl(2,\IC)\subset\sso(5,\IC)$ thus described is the one that is given by the Lefschetz operators
$L_\lambda$ and $\Lambda_\lambda$.
\silent{For the record $A\in \ssl_2$ defines $A\in\sso(3,\IC)$ by $C^3=\IC^2\oplus\IC e$,
$A(\alpha)=\alpha'\oplus \mu e$ with $\alpha\perp\alpha'$.
In terms of matrices it should be 
$$\xymatrix{
a&b\\
c&-a}
\to \xymatrix{0&b&a\\-b&0&c\\
-a&-c&0}$$}\qed
\bigskip

\bigskip

Let us conclude this section with a discussion of the moduli space of K3 surfaces without metrics and
generalized Calabi--Yau structures.
They are studied in terms of the following two period domains:
$$Q:=\{x~|~x^2=0,~x\bar x>0\}\subset \IP(H^2(M,\IC))$$
and
$$\widetilde Q:=\{x~|~\langle x,x\rangle=0,~\langle x,\bar x\rangle>0\}\subset \IP(H^*(M,\IC)),$$
where the latter involves the Mukai pairing $\langle~,~\rangle$. Clearly, $Q$ is naturally contained in
$\widetilde Q$.
\begin{definition}
We introduce the moduli space
$$\gn_{\rm gen}=\raisebox{.0ex}{$\{\varphi\cdot\IC\}$}_/\raisebox{-1.30ex}{$\cong$}$$ of all isomorphism classes
of generalized Calabi--Yau structures $\varphi\cdot \IC$ of hyperk\"ahler type on $M$ (cf.\ Definition \ref{defiso}). 
\end{definition}
By definition,  $\gn_{\rm gen}$ is the quotient of the set of  generalized Calabi--Yau structures of hyperk\"ahler type
on $M$ by the action of exact $B$-fields and of the groups $\IC^*$ and  ${\rm Diff}_*(M)$.
In contrast to the moduli spaces considered before, $\gn_{\rm gen}$ is not separated.

The classical counterpart is the moduli space 
$$\gn:=\raisebox{.0ex}{$\{\sigma
\cdot\IC\}$}_/\raisebox{-1.30ex}{${\rm Diff}_*(M)$}$$
of complex structures on $M$ or, equivalently, the moduli space of
marked K3 surfaces (or rather, one of the two connected components of it).

The classical period map $\kp:\gn\to Q$, $\sigma\cdot\IC\mapsto [\sigma]\cdot\IC$ extends
naturally  to $\kp_{\rm gen}:\gn_{\rm gen}\to \widetilde Q$,  $\varphi\cdot\IC\mapsto[\varphi]\cdot\IC$.
This yields a commutative diagram
$$
%\begin{array}{ccccccc}
\xymatrix@H=14pt{\gn\ar@{^(->}[d]\ar[r]^{\kp}& Q\ar@{^(->}[d]\ar@{^(->}[r]&\IP(H^2(M,\IC))\ar@{^(->}[d]\\
\gn_{\rm gen}\ar[r]^{\kp_{\rm gen}}&\widetilde Q\ar[r]&\IP(H^*(M,\IC)).}
%\end{array}
$$
The surjectivity of the period map
$\kp:\gn\to Q\subset\IP(H^2(M,\IC))$, a result due to Todorov, Looijenga, Siu (cf.\ \cite{Periodes}), is 
one of the fundamental results in the theory of K3 surfaces.
It easily generalizes to our situation:

\begin{proposition}
The period map $\kp_{\rm gen}$ is \'etale and
surjective. Moreover, $\kp_{\rm gen}$ is bijective over the
complement of the hyperplane section $\IP(H^2(M,\IC)\oplus H^4(M,\IC))\cap \widetilde Q$.
\end{proposition}

\prf The argument for the surjectivity follows the classical proof. We define an equivalence
relation on $\widetilde Q$ using positive four-spaces. Two positive planes $P,P'\in \widetilde Q$
are called equivalent if they generate a positive four-space. Obviously, this equivalence relation is open and,
since $\widetilde Q$ is connected, there exists only one equivalence class. In particular, it suffices
to show that for two planes $P,P'\in\widetilde Q$ generating a positive four-space $\Pi$ one has
$P\in{\rm Im}(\kp_{\rm gen})$ if and only if $P'\in{\rm Im}(\kp_{\rm gen})$. Since $\IP(\Pi_\IC)$ intersects $Q$, we may assume
that $P\in Q$ and that $\Pi$ is generated by $P$ and the plane $P_\varphi$ with $\varphi=\exp([B]+i[\omega])$
(cf.\ the proof of Prop.\ \ref{clgen}), where $[\omega]$ is a class in the positive cone. 
After changing $P$ a little in $Q$, which is allowed as the positivity of $\langle P,P'\rangle$
is preserved, we may assume that the Picard group of the corresponding
K3 surface is trivial. Hence, every element in the positive cone is actually a K\"ahler class.
This shows that the intersection $\IP(\Pi_\IC)\cap\widetilde Q$ is the image of the generalized twistor space and thus
contained in the image of $\kp$.

The last assertion follows from the Global Torelli theorem for generalized Calabi--Yau structures
$\varphi$ with $\varphi_0\ne0$ (see the proof of Proposition \ref{GT}).
\qed

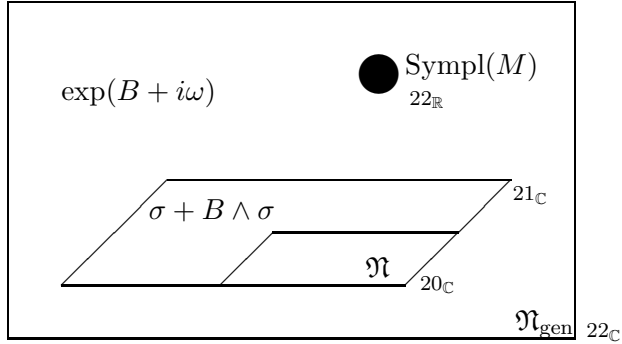
\bigskip

$$\put(-100,-10){\framebox(214,127)[br]{$\gn_{\rm gen}$}}
\put(115,-8){\makebox{${~}_{22_\IC}$}}
\put(40,90){\circle*{20}}
\put(48,80){\makebox{${~}_{22_\IR}$}}
\put(50,90){\makebox{${\rm Sympl}(M)$}}
\put(0,30){\line(5,-0){70}}
\put(-20,10){\line(5,-0){70}}
\put(-20,10){\line(1,1){20}}
%\put(50,10){\line(1,1){20}}
%\put(-20,5){\line(5,-0){90}}
\put(-80,10){\line(5,-0){120}}
\put(-80,10){\line(1,1){40}}
\put(50,10){\line(1,1){40}}
\put(-40,50){\line(5,-0){130}}
\put(35,13){\makebox{$\gn$}}
\put(52,10){\makebox{${~}_{20_\IC}$}}
\put(-80,80){\makebox{$\exp(B+i\omega)$}}
\put(87,44){\makebox{${~}_{21_\IC}$}}
\put(-47,35){\makebox{$\sigma+B\wedge\sigma$}}
$$

\bigskip

It is noteworthy that by incorporating  $B$-fields it is now possible to deform a K3 surface, i.e.\
a complex structure on $M$, continuously to a symplectic form.

One should think of $\IP(H^2(M,\IC)\oplus H^4(M,\IC))\cap \widetilde Q$ as the period domain
for the moduli space of all $B$-field shifts of $Q$. Thus, the latter has to be considered
as a hyperplane section of $\gn_{\rm gen}$. Note that the complement of this hyperplane section
is complex $22$-dimensional and parametrizes $B$-field shifts of hyperk\"ahler symplectic forms, which
for themselves form a real $22$-dimensional subspace.

%%%%%%%%%%%%%%%%%%%%%%%%%%%%%%%%%%%%%%%%%%%%%%%%%%%%%%%%%%%%%%%%%%%%%%%%%%%%%%%%%%%%%%%%%%%%%%
\section{Generalized Calabi--Yau structures and Derived Categories}\label{Orl}

The modest aim of the present section is to illustrate that 
generalized Calabi--Yau structures and generalized K3 surfaces seemingly provide a natu\-ral framework
for certain results and conjectures on equivalences of derived categories of coherent 
sheaves on algebraic K3 surfaces. More details can be found in \cite{HS}.

Let us begin with Orlov's result \cite{Orlov} on the equivalence of derived
cate\-gories of K3 surfaces. In the following 
let $X$ and $X'$ be two algebraic K3 surfaces given by Calabi--Yau structures
$\varphi=\sigma$ and $\varphi'=\sigma'$, respectively. Furthermore, let 
${\rm D}:={\rm D}^{\rm b}({\bf Coh}(X))$ and ${\rm D}':={\rm D}^{\rm b}({\bf Coh}(X'))$,
respectively, be their derived categories of coherent sheaves. 
Then Orlov's result can be reformulated as follows:

\begin{theorem}{\bf (Orlov)}
There exists an exact equivalence between the triangulated categories
${\rm D}$ and ${\rm D'}$  if and only if one of the following two
equivalent conditions is satisfied:

{\rm i)} There exists an Hodge isometry $T(X)\cong T(X')$ between their transcendental lattices.

{\rm ii)} The periods of $\varphi$ and $\varphi'$ are contained
in the same ${\rm O}(H^*(M,\IZ))$-orbit  of the natural action of the orthogonal
group of the Mukai lattice on the period domain.\qed
\end{theorem}

Orlov's result thus generalizes beautifully the Global Torelli theorem saying that
$X$ and $X'$ (not necessarily algebraic) are isomorphic if and only if their periods
are contained in the same ${\rm O}(H^2(M,\IZ))$-orbit. Thus, passing from
the period domain $Q$ to the generalized period domain $\widetilde Q$, or from 
${\rm O}(H^2(M,\IZ))$ to ${\rm O}(H^*(M,\IZ))$, corresponds in geometrical terms
to the passage from isomorphism classes of K3 surfaces to equivalence classes of their
derived categories. 

Note however that, for the time being, a completely satisfactory geome\-trical interpretation of
the action of ${\rm O}(H^*(M,\IZ))$ is missing. For the smaller orthogonal
group ${\rm O}(H^2(M,\IZ))$ it is provided by the surjectivity of the natural representation
${\rm Diff}(M)\to{\rm O}_+(H^2(M,\IZ))$ due to Borcea and Donaldson, where
${\rm O}_+\subset{\rm O}$ is the subgroup of orthogonal transformations that
preserve the orientation of the positive directions.

\medskip

What about generalized Calabi--Yau structures of the form $\exp(B)\sigma$? Here the situation
is less clear, but a conjecture treating this case has been proposed by A.\ C\u{a}ld\u{a}raru
in \cite{Cal}. 
Again, we assume that $X$ and $X'$ are two algebraic K3 surfaces. Moreover, we 
choose two torsion classes $\alpha\in H^2(X,\ko_X^*)$ and $\alpha'\in H^2(X',\ko_{X'}^*)$.
Lifting these classes to elements in $H^2(X,\IQ)$ respectively in 
$H^2(X',\IQ)$ one defines 
$T(X,\alpha):={\rm Ker}(\alpha:T(X)\to\IQ/\IZ)$ and similarly $T(X',\alpha')$.
(See below for more details.)

\begin{conjecture}{\bf (C\u{a}ld\u{a}raru)}
There exists an exact equivalence between the derived categories
${\rm D}^{\rm b}({\bf Coh}_\alpha(X))$ and ${\rm D}^{\rm b}({\bf Coh}_{\alpha'}(X'))$ of
twisted coherent sheaves if and only if there exists a Hodge isometry
$T(X,\alpha)\cong T(X',\alpha')$.
\end{conjecture}

This conjecture has been verified in a few special cases by C\u{a}ld\u{a}raru
himself \cite{Cal} and, more recently, by Donagi and Pantev \cite{DP}.
One way to define the derived category ${\rm D}^{\rm b}({\bf Coh}_\alpha(X))$
is to view $\alpha$ as an  Azumaya algebra $\ka$ (up to equivalence)
and to derive the abelian category of coherent $\ka$-modules. (That the Azumaya
algebra $\ka$ exists is due to Grothendieck, at least in the algebraic situation
considered here. The analytic analogue was recently established in 
\cite{HS}.)

In the rest of this section we indicate how to define the transcendental lattice
of a generalized Calabi--Yau structure and how to rephrase C\u{a}ld\u{a}raru's conjecture in terms
of periods and the ${\rm O}(H^*(M,\IZ))$-action.

\begin{definition}
The \emph{Picard group} of a generalized Calabi--Yau structure
$\varphi$ is the orthogonal complement of its
cohomology class:
$${\rm Pic}(\varphi):=\{\delta~|~\langle\delta,\varphi\rangle=0\}\subset H^*(M,\IZ).$$
\end{definition}

\begin{example}
i) If $\varphi$ is an ordinary Calabi--Yau structure $\sigma$, i.e.\
$\sigma$ is the holomorphic two form with respect to a complex structure
on $M$ defining a K3 surface $X$, then 
$${\rm Pic}(\varphi)=H^0(M,\IZ)\oplus {\rm Pic}(X)\oplus H^4(M,\IZ).$$
So, the Picard group of a K3 surface $X$ and the Picard group of the
Calabi--Yau structure naturally defined by it differ just by the hyperbolic
plane $H^0\oplus H^4$.

ii)  If $\omega$ is a symplectic structure and $\varphi=\exp(i\omega)$, then
$${\rm Pic}(\varphi)=H^2(M,\IZ)_\omega\oplus \{\delta_0+\delta_4\in (H^0\oplus H^4)(M,\IZ)~|~
\delta_0\int_M\omega^2=2\int_M\delta_4\}.$$
Here, $H^2(M,\IZ)_\omega$ is the group of $\omega$-primitive classes. Note that for $\omega$
very general the Picard group ${\rm Pic}(\exp(i\omega))$ is trivial, as $\int_M\omega^2$ will be irrational
and any integral class orthogonal to $\omega$ will be trivial.

iii) Let us twist an ordinary Calabi--Yau structure $\sigma$ by a B-field $B\in H^2(M,\IR)$,
i.e.\ we consider 
$\varphi:=\exp(B)\sigma=\sigma+ B^{0,2}\wedge\sigma$.
Then
$${\rm Pic}(\varphi)=H^4(M,\IZ)\oplus\{\delta_0+\delta_2\in (H^0\oplus H^2)(M,\IZ)
~|~\int_M\delta_2\wedge\sigma=\delta_0\int_MB\wedge \sigma\}.$$
Clearly, $H^4(M,\IZ)\oplus {\rm Pic}(X)\subset {\rm Pic}(\exp(B)\sigma)$, where
$X$ is the K3 surface defined by $\sigma$. In general, this inclusion will be strict.
Note that ${\rm Pic}(\varphi)$ depends only on $\sigma$ and the
$(0,2)$-part of $B$.
\end{example}

If $X$ is a classical K3 surface then the Hodge decomposition $H^2(X,\IC)=H^{2,0}(X)
\oplus H^{1,1}(X)\oplus H^{0,2}(X)$ induces a weight two Hodge structure of
$H^*(M,\IZ)$ whose $(1,1)$-part is $H^0(X,\IC)\oplus H^{1,1}(X)\oplus H^4(X,\IC)$
and its $(2,0)$-part is spanned by the Calabi--Yau form $\sigma$.

In the same vain, any generalized Calabi--Yau structure $\varphi$ defines a
weight two Hodge structure on 
$H^*(M,\IZ)$ whose $(2,0)$-part is spanned by the cohomology class
of $\varphi$. This determines the other parts by requiring that the Hodge
decomposition is orthogonal with respect to the Mukai pairing.
Clearly, the Picard group ${\rm Pic}(\varphi)$ 
is of pure type $(1,1)$.

\begin{definition}
The \emph{transcendental lattice} $T(\varphi)$ of a generalized Calabi--Yau structure
 $\varphi$ on $M$
is the orthogonal complement of ${\rm Pic}(\varphi)$ in $H^*(M,\IZ)$.
More precisely,
$$T(\varphi):=\{\gamma\in H^*(M,\IZ)~|~\langle\gamma,\delta\rangle=0~{\rm for~all~}\delta\in
{\rm Pic}(\varphi)\}.$$
\end{definition}

Via the intersection pairing, the transcendental lattice $T(\varphi)$ can be identified
with the dual of $H^*(M,\IZ)/{\rm Pic}(\varphi)$. Also note that the transcendental lattice
$T(\varphi)$ is almost never pure.

\begin{example}
i) Let us recall that the transcendental lattice $T(X)$ of a classical K3 surface $X$
is the orthogonal complement of ${\rm Pic}(X)$ inside $H^2(X,\IZ)$.
One easily sees that $T(X)=T(\sigma)$ for a Calabi--Yau form $\sigma$ on
$X$.

iii) For $\varphi=\exp(B)\sigma$ we find that $\gamma_0+\gamma_2+\gamma_4\in
T(\varphi)$ implies $\gamma_0=0$ and $\gamma_2\in T(X)$.
\end{example}

Now we shall explain how to identify the transcendental lattice
$T(\varphi)$ in the case $\varphi=\exp(B)\sigma$ with $T(X,\alpha_B)$ defined before.
Here, $B\in H^2(X,\IR)$ is a B-field whose $(0,2)$-component $B^{0,2}\in H^{0,2}(X)=H^2(X,\ko)$
induces a torsion element $\alpha_B$ in $H^2(X,\ko^*)$ via the exponential
map $H^2(X,\ko)\to H^2(X,\ko^*)$. As before, $X$ is the K3 surface defined by $\sigma$.

Let us first recall some basic facts concerning the B-field and its associated 
Brauer class $\alpha_B$:

i) Since $B$ is real,  $B^{0,2}\in H^2(X,\ko)$ is trivial if and only if $B\in H^{1,1}(X,\IR)$.

ii)  The $(0,2)$-part  of $B$ is contained in the image of $H^2(X,\IZ)\to H^2(X,\ko)$
if and only if $B\in H^2(X,\IZ)+H^{1,1}(X,\IR)$.

iii) Using the exponential sequence $H^2(X,\IZ)\to H^2(X,\ko)\to H^2(X,\ko^*)$,
one finds that the B-field $B\in H^2(X,\IR)$ defines an $r$-torsion class
$\alpha_B\in H^2(X,\ko^*)$ if and only if $rB\in H^2(X,\IZ)+H^{1,1}(X,\IR)$.

Let us now assume that $\alpha_B\in H^2(X,\ko^*)$ is $r$-torsion. Then we can
write $B=\beta+\delta$ with $\beta\in (1/r)H^2(X,\IZ)$ and $\delta\in H^{1,1}(X,\IR)$.
This decomposition is not unique. If, however, $B\in H^2(M,\IQ)$ then
$\delta\in{\rm Pic}(X)_\IQ$ and, hence, the induced homomorphism
$B:T(X)\to \IQ$, $\gamma\mapsto \int_M\gamma\wedge B$ can be described
in terms of $\beta$ as $\gamma\mapsto\int_M\gamma\wedge\beta$. Clearly,
the image is contained $(1/r)\IZ\subset\IQ$ and one has
$$T(X,\alpha_B)={\rm Ker}(B:T(X)\to\IQ/\IZ)=\{\gamma\in T(X)~|~\gamma\beta\in\IZ\}.$$
Note that if $B$ is not rational, even if the induced $\alpha_B$ is torsion, then the
map $\gamma\mapsto\int_M\gamma\wedge B$ takes values in $\IR$ and usually $T(X,\alpha_B)$ defined
as the kernel of the induced map to $\IR/\IZ$ will be too small to be interesting.
But if one starts out with a torsion class $\alpha\in H^2(X,\ko^*)$ one always
finds a rational B-field $B$ with $\alpha=\alpha_B$. So,
we will continue to assume that $B$ is rational.

The sublattice $T(X,\alpha_B)\subset T(X)$ will be viewed as a sublattice of $H^*(M,\IZ)$
via the injection
$$\eta:T(X,\alpha_B)\to (H^2\oplus H^4)(M,\IZ), ~~\gamma\mapsto \gamma+\gamma\wedge B.$$
Note that $\eta$ is the restriction of the isometry (with respect to the Mukai pairing)
$$\xymatrix{\exp(B):H^*(X,\IQ)\ar[r]^-\sim&H^*(X,\IQ)}.$$

\begin{proposition}
For a rational B-field $B$ the map $\eta$ defines
an Hodge iso\-metry 
$$T(X,\alpha_B)\cong T(\varphi),$$
where $\varphi$ is the generalized Calabi--Yau structure $\sigma+B\wedge \sigma$.
\end{proposition}

\begin{proof}
Firstly, it is clear that $\eta$ is compatible with the quadratic forms provided
by the intersection and the Mukai pairing, respectively.

Secondly, the Hodge structures are respected. Indeed, the $(2,0)$-part
of $T(X,\alpha_B)$ is generated by $\sigma$ and $\eta(\sigma)=\varphi$ spans
the $(2,0)$-part of the Hodge structure on $H^*(M,\IZ)$ defined by 
$\varphi$. (In fact, $\eta$ is nothing but $\exp(B)$, so it is clearly compatible
with the Hodge structures.)

Thus, it suffices to show that $\eta$ is indeed bijective. Let us first verify the inclusion
$T(\varphi)\subset \eta(T(X,\alpha_B))$:
Let $\gamma=\gamma_0+\gamma_2+\gamma_4\in T(\varphi)$. As was observed above,
one has $\gamma_0=0$ and $\gamma_2\in T(X)$. As $T(X,\alpha_B)$ and $T(\varphi)$ are independent
of the $(1,1)$-part of $B$, we may assume $B=\beta\in(1/r)H^2(M,\IZ)$.
Then $r+rB+0\in{\rm Pic}(\varphi)$, for $\int_M(rB)\wedge\sigma=r\int_M B\wedge\sigma$.
Thus, $\gamma$ is orthogonal to $r+rB$ and, therefore,
$\int_M\gamma_2\wedge (rB)=r\int_M\gamma_4$. Hence, $\gamma\in{\rm Im}(\eta)$.

For the other inclusion one has to check that $\gamma:=\gamma_2+\gamma_2\wedge B$
is orthogonal to ${\rm Pic}(\varphi)$ for any $\gamma_2\in T(X,\alpha_B)$.
Clearly, any such $\gamma$ is orthogonal to $H^4(M,\IZ)\oplus {\rm Pic}(X)$
and to the special element $r+rB\in{\rm Pic}(\varphi)$.
Let $\delta_0+\delta_2\in {\rm Pic}(\varphi)$. Then $\delta_0\int_MB\wedge\sigma=
\int_M\delta_2\wedge\sigma$. Writing $B^{0,2}=\mu\bar\sigma$
and $\delta_2^{0,2}=\lambda\bar\sigma$, this condition becomes
$\delta_0\mu=\lambda$. Hence, $\delta_2=\delta_0B+\delta_2'$ with $\delta_2'\in
H^{1,1}(X,\IQ)$. Here, the rationality of $\delta_2'$ follows from the
rationality of the other terms.
For $\gamma_2\in T(X)$ this yields $\gamma_2\wedge\delta_2=\delta_0\gamma_2\wedge B+
\gamma_2\wedge\delta_2'=\delta_0\gamma_2\wedge B$, i.e.\ $\eta(\gamma_2)$ is
orthogonal to $\delta_0+\delta_2$.
\end{proof}

\begin{corollary}
There exists an Hodge isometry $T(X,\alpha_B)\cong T(X',\alpha_{B'})$
if the periods of the two generalized Calabi--Yau structures
$\varphi=\sigma+B\wedge\sigma$ and $\varphi'=\sigma'+B'\wedge \sigma'$ are contained
in the same orbit of the natural action of ${\rm O}(H^*(M,\IZ))$. 
\end{corollary}

\begin{proof}
This is obvious, as $T(X,\alpha_B)$ and $T(\varphi)$ are Hodge isometric and
$T(\varphi)$ is defined
in terms of the period.
\end{proof}

In many cases the corollary can be improved to an ``if and only if'' statement. In general,
however, this is not true. In the untwisted case any Hodge isometry
$T(X)\cong T(X')$ can be extended to an Hodge isometry $\widetilde H(X,\IZ)\cong
\widetilde H(X',\IZ)$ due to results of Nikulin \cite{Nik}. The existence of the hyperbolic
plane $H^0\oplus H^4$ in the orthogonal complement $T^\perp$ is crucial for this. In the twisted
case, $T(\varphi)^\perp$ does not necessarily contain such an hyperbolic plane. In fact,
in \cite{HS} we give an explicit example of an Hodge isometry $T(\varphi)\cong T(\varphi')$
which cannot be extended simply due to the fact that the Picard groups
are not isometric. We furthermore put forward a version of C\u{a}ld\u{a}raru's conjecture
that takes into account not only the transcendental lattice
$T(X,\alpha_B)\cong T(\varphi)$, but the full Hodge structure defined by the generalized
Calabi--Yau structure $\varphi=\sigma+B\wedge\sigma$.

\bigskip

Clearly, the action of the discrete group ${\rm O}(H^*(M,\IZ))$ studied in this section
is intimately related to mirror symmetry for K3 surfaces. For instance one can observe
that ${\rm O}(H^*(M,\IZ))$ frequently interchanges the periods of honest K3 surfaces
with those of symplectic structures. Although mirror symmetry phenomena on the level
of moduli spaces are much simpler for K3 surfaces than for arbitrary Calabi--Yau manifolds, not
much is known about it on a deeper level, e.g.\  mirror symmetry in its homological version
(but see the recent paper \cite{Seidel}) or mirror symmetry as a duality for conformal
field theories, vertex algebras as in \cite{KO} for tori.

%%%%%%%%%%%%%%%%%%%%%%%%%%%%%%%%%%%%%%%%%%%%%%%%%%%%%%%%%%%%%%%%%%%%%%%%%%%%%%%%%%%%%%%%%%%%%%%%%
\section{The moduli space of symplectic structures as a Lagrangian}

This section is devoted to a canonical symplectic form on the moduli space $\gn_{\rm gen}$ of generalized
Calabi--Yau structures on $M$. We will show that it coincides with the pull-back (via the period 
map) of the curvature of the tautological bundle with respect to a certain hermitian structure. Moreover, the
subset of generalized Calabi--Yau structures of the form $\exp(i\omega)$ is Lagrangian.

As before, $\langle~,~\rangle$ denotes the Mukai pairing on even forms.

\begin{lemma}
The pairing $H(\varphi,\psi)=\int\langle\varphi,\overline\psi\rangle$ defines a non-degenerate indefinite hermitian
product on $\ka_\IC^{2*}(M)$.
\end{lemma}

\prf E.g.\ if $\varphi_2\ne0$, then $\int\langle\varphi,*\varphi_2\rangle=||\varphi_2||^2>0$.
Together with similar calculations in the other cases this shows that $H$ is non-degenerate.
Moreover, $H$ is indefinite, for $H(\omega,\omega)>0$ for any symplectic structure $\omega$
and $H(\alpha,\alpha)<0$ for any primitive $(1,1)$-form
$\alpha$ (with respect to an arbitrary Calabi--Yau structure and a K\"ahler form).
\qed

\begin{corollary}
The constant two-form $\Omega:={\rm Im}(H)$ defines a symplectic structure on $\ka_\IC^{2*}(M)$.\qed
\end{corollary}

Since $\int\langle~,~\rangle$ is invariant under $B$-field transformations and diffeomorphisms,
we obtain this way a hermitian structure $H$  and a two-form $\Omega$ on $\gn_{\rm gen}$.
That $\Omega$ is indeed a symplectic structure can be seen by the following explicit description.

\begin{lemma}
Let $\varphi$ be a generalized Calabi--Yau structure on $M$. Then the tangent space
$T_\varphi\gn_{\rm gen}$ can be naturally identified with
$$\{\alpha\in H^*(M,\IC)~|~\langle\alpha,[\varphi]\rangle=\langle\alpha,[\overline\varphi]\rangle=0\}.$$
The form $\Omega$, which is given by $\Omega(\alpha,\beta)={\rm Im}\langle\alpha,\bar\beta\rangle$,
defines a symplectic structure, i.e.\ it is  non-degenerate and closed.
\end{lemma} 
\prf
One way to prove this, is to note that the period map $\kp_{\rm gen}$ is a local isomorphism (this is the
analogue of the Local Torelli theorem). Hence, it suffices to compute the tangent space of $\widetilde Q$ at
$[\varphi]$. Now, $\langle [\varphi]+\varepsilon \alpha,[\varphi]+\varepsilon \alpha\rangle=0$ yields
$\langle [\varphi],\alpha\rangle=0$. Thus, the tangent space $T_{[\varphi]}\widetilde Q$ is canonically
isomorphic to $[\varphi]^\perp/\IC\cdot[\varphi]$, which can be identified with the subspace given above, due 
to $\langle[\varphi],[\overline\varphi]\rangle\ne0$.

The Mukai pairing $\langle~,~\rangle$ is non-degenerate on the orthogonal complement of $P_\varphi$, i.e.\
on the tangent space $T_\varphi\gn_{\rm gen}$. Hence, $\Omega$ is a non-degenerate two-form on $\gn_{\rm gen}$.
 
The closedness of $\Omega$ follows from the construction: Its pull-back is a restriction of a constant form
on $\ka^{2*}(M)_\IC$.
\qed

\bigskip

Note that since $H$ is indefinite, the symplectic structure $\Omega$ is not K\"ahler with respect to the
natural complex structure on $\gn_{\rm gen}$. 

Let $\ko(-1)$ be the tautological line bundle on $\IP(H^*(M,\IC))$. The hermi\-tian structure on $H^*(M,\IC)$
defined by the Mukai pairing $\langle\alpha,\bar\beta\rangle$ can be viewed as a constant hermitian structure
on the constant bundle $H^*(M,\IC)\otimes\ko$ on $\IP(H^*(M,\IC))$. Consider the Euler sequence and, in particular,
the natural inclusion $\ko(-1)\subset H^*(M,\IC)\otimes\ko$, which provides   a natural hermitian structure on
$\ko(-1)$. Then the constant connection on $H^*(M,\IC)\otimes\ko $
induces the Chern connection $\nabla$ on $\ko(-1)$. A standard calculation, well-known
in the positive definite case, relates the curvature
$F_\nabla$ to the above defined symplectic form $\Omega$:

\begin{proposition}
$\kp^*(iF_\nabla)=-\Omega$.\qed
\end{proposition}
\silent{This is analogous  to one what one does for the Fubini-Study metric.
One endows a vector space $V$ with a hermitian product. This yields a K\"ahler structure on $\IP(V)$. The
curvature of $\ko(1)$ is nothing but the Fubini-Study K\"ahler form. On $\ko(-1)$ one acquires the sign.\qed

\bigskip}

As remarked above, we cannot expect, due to the indefiniteness of the chosen hermitian structure on $\IP(H^*(M,\IC))$,
that the curvature satisfies any positivity condition.

\bigskip

Let ${\rm Sympl}(M)$ be the moduli space of symplectic structures on $M$, which is identified with
a submanifold of $\gn_{\rm gen}$ via the 
the natural inclusion $\omega\mapsto\exp(i\omega)$.
(For simplicity we assume here that any symplectic structure on $M$ is K\"ahler.
Otherwise, one has to work with the moduli space of
all generalized Calabi--Yau structures and not only of those of hyperk\"ahler type.)
An easy dimension count shows that ${\rm Sympl}(M)$ is
a real $22$-dimensional submanifold of the complex $22$-dimensional complex manifold $\gn_{\rm gen}$.
\begin{proposition}
The submanifold ${\rm Sympl}(M)\subset \gn_{\rm gen}$ is a Lagrangian submanifold with respect
to $\Omega$.
\end{proposition}
\prf
The tangent space $T_\omega{\rm Sympl}(M)$ is canonically identified with
$H^2(M,\IR)$ and the tangent map $d\exp:T_{\omega}{\rm Sympl}(M)\to T_{\exp(i\omega)}$ is 
given by
$$\begin{array}{ccl}
H^2(M,\IR)&\rpfeil{5}{}& T_{\exp(i\omega)}\gn_{\rm gen}%\\
%&&=\{\alpha~|~\langle\alpha,[\exp(i\omega)]\rangle=\langle\alpha,[\exp(-i\omega)]\rangle=0\}
\subset H^*(M,\IC)\\
\alpha&\mapsto &i\alpha-\alpha\omega\\
\end{array}$$
Clearly, $\Omega(i\alpha-\alpha\omega,i\beta-\beta\omega)={\rm Im}\langle i\alpha,-i\beta\rangle=0$,
whenever $\alpha,\beta$ are both real cohomology classes.
\qed

\begin{remark}
By construction the curvature $F_\nabla$ is invariant under the action of the orthogonal group
$\OO(H^*(M,\IR))=\OO(4,20)$. Thus, translating ${\rm Sympl}(M)$ by elements in $\OO(4,20)$ yields new
Lagrangian submanifolds.
\silent{Are there other interesting ones? Is ${\rm Sympl}$ in any way distinguished?}
\end{remark}

\bigskip

\bigskip


\vskip1cm

{\footnotesize \begin{thebibliography}{mm}

\bibitem{Periodes} \em G\'eom\'etrie des surfaces K3: modules et p\'eriodes.
\em S\'eminaires Palaiseau. ed A.\ Beauville, J.-P.\ Bourguignon, M.\ Demazure.
Ast\'erisque 126 (1985).

\bibitem{AM} P.\ Aspinwall, D.\ Morrison \em String theory on $K3$ surfaces. \em  Mirror symmetry II AMS/IP
Stud.\ Adv.\ Math.\ 1 (1997), 703-716.
 hep-th/9404151.
 
\bibitem{Borcea} C.\ Borcea \em Diffeomorphisms of a K3 surface. \em
Math.\ Ann.\ 275 (1986), 1-4.

\bibitem{Cal} A.\ C\u{a}ld\u{a}raru
\em Derived categories of twisted sheaves on Calabi--Yau manifolds. \em
Ph.-D. thesis Cornell (2000).

\bibitem{DP} R.\ Donagi, T.\ Pantev
\em Torus fibrations, gerbes, and duality. \em
math.AG/0306213

\bibitem{Donaldson} S.\ Donaldson \em Moment maps and diffeomorphisms. \em Asian J.\ Math.\ 3 (1999), 1-16.

\bibitem{Donalson2} S.\ Donaldson \em Polynomial invariants for smooth four-manifolds. \em
 Top.\ 29 (1990), 257-315.

\bibitem{FM} R.\ Friedman, J.\ Morgan \em Smooth four-manifolds and complex surfaces. \em
Erg.\ Math.\ 27 (1994), Springer.

\bibitem{Gualtieri} M.\ Gualtieri
\em Generalized Complex Geometry. \em math.DG/0401221.Ph.D.-thesis Oxford.


\bibitem{Hitchin} N.\ Hitchin \em Generalized Calabi--Yau manifolds. \em 
 Q.\ J.\ Math.\  54 (2003), 281-308.

\bibitem{Huy} D.\ Huybrechts \em Moduli spaces of hyperk\"ahler manifolds and mirror symmetry. \em
School on Intersection Theory and Moduli. Trieste 2002.

\bibitem{HS} D.\ Huybrechts, S.\ Schr\"oer
\em The Brauer group of analytic K3 surfaces. \em
IMRN 50 (2003), 2687-2698.

\bibitem{HS} D.\ Huybrechts, P.\ Stellari
\em Equivalences of twisted K3 surfaces. \em
math.AG/0409030.

\bibitem{KO}  A.\ Kapustin, D.\ Orlov
\em Vertex algebras, mirror symmetry, and D-branes: the case of complex tori. \em
Comm.\ Math.\ Phys.\  233 (2003), 79-136.

\bibitem{McDuffSal} D.\ McDuff, D.\ Salamon
\em Introduction to symplectic topology. \em 
Oxford University Press, (1998). 

\bibitem{Moser} J.\ Moser \em On the volume elements on a manifold. \em Trans.\ AMS 120 (1965), 286-294.

\bibitem{NW} W.\ Nahm, K.\ Wendland \em A hiker's guide to $K3$. Aspects of $N=(4,4)$ superconformal
field theory with central charge $c=6$. \em
Comm.\ Math.\ Phys.\ 216 (2001), 85-138.

\bibitem{Nik} V.V.\ Nikulin
\em Integral symmetric bilinear forms and some of their
applications. \em
 Math USSR Izvestija 14 (1980), 103-167.

\bibitem{Orlov} D.\ Orlov, \em Equivalence of derived categories and
K3 surfaces. \em Alg.\ Geom.\ 7, J.\ Math.\ Sci.\ (New York)
84 (1997), 1361-1381. 

%\bibitem{Seidel} P.\ Seidel \em Floer homology and the symplectic isotopy problem. \em Ph.D. thesis
%Oxford

\bibitem{Seidel} P.\ Seidel
\em Homological mirror symmetry for the quartic surface. \em
math.SG/0310414

\bibitem{Siu} Y.-T.\ Siu \em Every K3 surface is K\"ahler. \em Invent.\ math. 73 (1983), 139-150.

\bibitem{Verb} M.\  Verbitsky \em Cohomology of compact hyper-K\"{a}hler manifolds and its applications. \em 
 Geom.\ Funct.\ Anal.\ 6  (1996), 601-611. 
\bibitem{Yau} S.-T.\ Yau \em On the Ricci curvature of a compact K\"ahler manifold and the complex Monge-Amp\`ere
equation. \em Comm.\ Pure Applied Math.\ 31 (1978), 339-411.

\end{thebibliography}}
\end{document}